\documentclass[preprint,10.5pt]{elsarticle}

\usepackage{amssymb}
\usepackage{amsthm}

\usepackage{multicol}
\usepackage{mathrsfs}
\usepackage{amsfonts}
\usepackage{CJK,fancybox,color,graphicx,amsmath,amsthm,amssymb,enumerate}

\journal{}

\newtheorem{thm}{Theorem}[section]

\newtheorem{lem}[]{Lemma}[section]

\theoremstyle{definition}
\newtheorem{defn}[thm]{Definition}
\theoremstyle{remark}
\newtheorem{rem}[]{Remark}[section]
\numberwithin{equation}{section}
\newcommand{\D}{\displaystyle}
\newcommand{\DF}[2]{\frac{\D#1}{\D#2}}

\begin{document}

\begin{frontmatter}



\title{\bf Stability of Riemann solutions to pressureless Euler equations with Coulomb-like friction by flux approximation}

\author{Qingling Zhang}

\address{School of Mathematics and Computer Sciences, Jianghan University, Wuhan 430056, PR China }
\ead{zhangqingling2002@163.com}
\begin{abstract}

We study the stability of Riemann solutions to pressureless Euler
equations with Coulomb-like friction under the nonlinear
approximation of flux functions with one parameter. The approximated
system can be seen as the generalized Chaplygin pressure Aw-Rascle
model with Coulomb-like friction, which is also equivalent to the
nonsymmetric system of Keyfitz-Kranzer type with generalized
Chaplygin pressure and Coulomb-like friction. Compared with the
original system, The approximated system is strictly hyperbolic,
which has one eigenvalue genuinely nonlinear and the other linearly
degenerate. Hence, the structure of its Riemann solutions is much
different from the ones of the original system. However, it is
proven that the Riemann solutions for the approximated system
converge to the corresponding ones to the original system as the
perturbation parameter tends to zero, which shows that the Riemann
solutions to nonhomogeneous pressureless Euler equations is stable
under such kind of flux approximation. The result in this paper
generalizes the stability of Riemann solutions with respect to flux
perturbation from the well-known homogeneous case to the
nonhomogeneous case.

\end{abstract}


\end{frontmatter}

% \linenumbers


\section{Introduction}
\setcounter{equation}{0}

Non-strictly hyperbolic system  have important physical background,
which is also difficult and interesting in mathematics and attract
many people to study them. It is well known that their Cauchy
problem usually does not have a weak $L^{\infty}$-solution. A
typical example is the Cauchy problem for pressureless Euler
equations (which is also called as zero pressure flow or
transportation equations) \cite{Huang-Wang,Wang-Huang-Ding}.
Therefore, the measure-value solution should be introduced to this
nonclassical situation, such as delta shock wave
\cite{Chen-Liu1,Sheng-Zhang,Tan-Zhang-Zheng} and singular shock
\cite{Keyfitz-Kranzer,Nedeljkov}, which can also provide a
reasonable explanation for some physical phenomena. However, the
mechanism for the formation of delta shock wave cannot be fully
understood, although the necessity of delta shock wave is obvious
for Riemann solutions to some non-strictly hyperbolic system. Now
there are some related results for homogenous equations
\cite{Chen-Liu1,Shen-Sun}, but few results have been shown for
nonhomogeneous equations.

 \ \  In this paper, we are mainly
concerned with zero pressure flow with Coulomb-like friction
\begin{equation}\label{1.1}
\left\{\begin{array}{ll}
\rho_t+(\rho u)_x=0,\\
(\rho u)_t+(\rho u^{2})_x=\beta\rho,
 \end{array}\right.
\end{equation}
where the state variable $\rho>0$, $u$ denote the density and
velocity, respectively, and $\beta$ is a frictional constant.

 The motivation of study (\ref{1.1}) comes from the violent
 discontinuities in shallow flows with large Froude number \cite{Edwards-etc.}. It
 can also be derived directly from the so-called pressureless
 Euler/Euler-Possion systems \cite{Nguyen-Tudorascu}. Moreover, the system (\ref{1.1}) can
 also be obtained formally from the model proposed by Brenier et al.\cite{Brenier-etc.} to
 describe the sticky particle dynamics with interactions. Recently,
 the Riemann problem and shadow wave for (\ref{1.1}) have been studied respectively in
 \cite{Shen1}
 and \cite{Daw-Nedeljkov}. Remarkably, in \cite{Shen1}, it is shown that the Riemann problem for the nonhomogeneous
equations (\ref{1.1}) has delta shock wave solutions in some
situations.

Delta shock wave is a kind of nonclassical nonlinear wave on which
at least one of the state variables becomes a singular measure.
Korchinski \cite{Korchinski} firstly introduced the concept of the
$\delta$-function into the classical weak solution in his
unpublished Ph.D. thesis. In 1994, Tan, Zhang and Zheng
\cite{Tan-Zhang-Zheng} considered some 1-D reduced system and
discovered that the form of $\delta$-functions supported on shocks
was used as parts in their Riemann solutions for certain initial
data. Since then, delta shock wave has been widely investigated, see
\cite{Bouchut,Li-Zhang-Yang,Sheng-Zhang} and references cited
therein.

The formation of delta shock wave has been extensively studied by
the vanishing pressure approximation for zero pressure flow
\cite{Chen-Liu1,Shen-Sun} and Chaplygin gas dynamics
\cite{Cheng-Yang,Sheng-Wang-Yin,Yang-Wang1}. Recently, the flux
approximation with two parameters \cite{Yang-Wang2} and three
parameters \cite{Yang-Liu} has also been carried out for zero
pressure flow. In the present paper, we consider the nonlinear
approximation of flux functions for zero pressure flow with
coulomb-like friction which has not been paid attention before.

Specifically, we introduce the nonlinear approximation of flux
functions in (\ref{1.1}) as follows:
\begin{equation}\label{1.2}
\left\{\begin{array}{ll}
\rho_t+(\rho u)_x=0,\\
(\rho( u+P))_t+(\rho u( u+P))_x=\beta\rho,
 \end{array}\right.
\end{equation}
where $P$ is given by the state equation for generalized Chaplygin
gas \cite{Bilic-Tupper-Viollier,Setare,Sheng-Wang-Yin,Wang}
\begin{equation}\label{1.3}
P=-\frac{A}{\rho^{\alpha}}, \ \ A>0, \ 0<\alpha<1,
\end{equation}
with $\alpha$ a real constant and the parameter $A$ sufficiently
small. System (\ref{1.2}) and (\ref{1.3}) can be seen as the
generalized Chaplygin pressure Aw-Rascle model with Coulomb-like
friction. By taking $u=w-P$, (\ref{1.2}) can be written
 as follows:
\begin{equation}\label{1.4}
\left\{\begin{array}{ll}
\rho_t+(\rho (w-P))_x=0,\\
(\rho w)_t+(\rho w( w-P))_x=\beta\rho,
 \end{array}\right.
\end{equation}
with a pure flux approximation. (\ref{1.4}) together with
(\ref{1.3}) can also be seen as the nonsymmetric system of
Keyfitz-Kranzer type with generalized Chaplygin pressure and
Coulomb-like friction \cite{Guo}. Recently, for $\beta=0$, Cheng has
shown that the structure of the Riemann solutions to (\ref{1.2}) and
(\ref{1.4}) were very similar \cite{Cheng1,Cheng2}.

More precisely, we are only concerned with the Riemann problem, i.e.
the initial data taken as follows:
\begin{equation}\label{1.5}
(\rho,u)(x,0)=\left\{\begin{array}{ll}
(\rho_-,u_-),\ \ x<0,\\
(\rho_ +,u_+),\ \ x>0,
\end{array}\right.
\end{equation}
where $\rho_{\pm}$ and $u_{\pm}$ are all given constants.

 In this paper, we will find that the delta shock wave also appears
in the Riemann solutions to (\ref{1.2}) for some specific initial
data. We are interested in how the delta-shock solution of
(\ref{1.2}) and (\ref{1.5}) develops under the influence of the
Coulomb-like friction. The advantage of this kind source term is in
that (\ref{1.2}) can be written in a conservative form such that
exact solutions to the Riemann problem (\ref{1.2}) and (\ref{1.5})
can be constructed explicitly. We shall see that the Riemann
solutions to (\ref{1.2}) and (\ref{1.5}) are not self-similar any
more, in which the state variable $u$ varies linearly along with the
time $t$ under the influence of the Coulomb-like friction. In other
words, the state variable $u-\beta t$ remains unchanged in the left,
intermediate and right states. In some situations, the delta-shock
wave appears in the Riemann solutions to (\ref{1.2}) and
(\ref{1.5}). In order to describe the delta-shock wave, the
generalized Rankine-Hugoniot conditions are derived and the exact
position, propagation speed and strength of the delta shock wave are
obtained completely. It is shown that the Coulomb-like friction term
make contact discontinuities, shock waves, rarefaction waves and
delta shock waves bend into parabolic shapes for the Riemann
solutions.

Furthermore, it is proven rigorously that the limits of Riemann
solutions to (\ref{1.2}) and (\ref{1.5}) converge to the
corresponding ones to (\ref{1.1}) and (\ref{1.5}) when the
perturbation parameter $A$ tends to zero. In other words, the
Riemann solutions (\ref{1.1}) and (\ref{1.5}) is stable with respect
to the nonlinear approximations of flux functions in the form of
(\ref{1.2}). Actually, for the case $\alpha=1$ in (\ref{1.3}),
system (\ref{1.2}) becomes the Chaplygin pressure Aw-Rascle model
with Coulomb-like friction \cite{Pan-Han}. Similar result can be
easily got, so we do not focus on it here. Moreover, the results got
in this paper can also be generalized to the nonsymmetric system of
Keyfitz-Kranzer type (\ref{1.4}) with the same generalized Chaplygin
pressure and Coulomb-like friction.

This paper is organized as follows. In section 2, we describe simply
the solutions of the Riemann problem (\ref{1.1}) and (\ref{1.5}) for
completeness. In Section 3, the approximated system (\ref{1.2}) is
reformulated into a conservative form and some general properties of
the conservative form are obtained. Then, the exact solution to the
Riemann problem for the conservative form are constructed
explicitly, which involves the delta shock wave. Furthermore, the
generalized Rankine-Hugoniot conditions are established and the
exact position, propagation speed and strength of the delta shock
wave are given explicitly. In Section 4, the generalized
Rankine-Hugoniot conditions and three kinds of Riemann solutions to
the approximated system (\ref{1.2}) and (\ref{1.5}) are given.
Furthermore, it is proven rigorously that the delta-shock wave is
indeed a week solution to the Riemann problem (\ref{1.2}) and
(\ref{1.5}) in the sense of distributions. In Section 5, the limit
of Riemann solutions to the approximated system (\ref{1.2}) is taken
by letting the perturbation parameter $A$ tends to zero, which is
identical with the corresponding ones to the original system.
Finally, conclusions and discussions are drawn in Section 6.

\section{preliminaries}
In this section, we simply describe the results on the Riemann
problem (\ref{1.1}) and (\ref{1.5}), which can be referred to
\cite{Shen1} in details.

Let us first state some known fact about elementary waves of the
given system. The system (\ref{1.1}) is weakly hyperbolic with the
double eigenvalue $\lambda_{1}=\lambda_{2}=u$. Let us first look for
a solution to (\ref{1.1}) when initial data are constants,
$(\rho(x,0),u(x,0))=(\rho_{0},u_{0})$. For smooth solutions, one can
substitute $\rho_{t}$ from the first equation of (\ref{1.1}) into
the second one and eliminate $\rho$ from it by division (provieded
that we are away from a vacuum state). So, we have now the equation
$u_{t}+uu_{x}=\beta$ that can be solved by the method of
characteristics: $u=u_{0}+\beta t, x=x_{0}+u_{0}t+\frac{1}{2}\beta
t^{2}$. The first equation then becomes $\rho_t+(u_{0}+\beta
t)\rho_{x}=0$ with a solution $\rho=\rho_{0}$ on each curve $
x=x_{0}+u_{0}t+\frac{1}{2}\beta t^{2}$. So, the solution for
constant initial data is $(\rho,u)=(\rho_{0},u_{0}+\beta t)$.

For the case $u_-<u_+$, there is no characteristic passing through
the region $\{(x,t): u_-t+\frac{1}{2}\beta
t^{2}<x<u_+t+\frac{1}{2}\beta t^{2}\}$, so the vacuum should appear
in the region. The solution can be expressed as
\begin{equation}\label{2.1}
(\rho,u)(x,t)=\left\{\begin{array}{ll}
(\rho_-,u_-+\beta t),\ \ \ \ -\infty<x<u_-t+\frac{1}{2}\beta t^{2},\\
vacuum,\ \ \ \ \ \ \ \ \ \ u_-t+\frac{1}{2}\beta t^{2}<x<u_+t+\frac{1}{2}\beta t^{2},\\
(\rho_+,u_++\beta t),\ \ \ \ u_+t+\frac{1}{2}\beta t^{2}<x<\infty.
\end{array}\right.
\end{equation}

For the case $u_-=u_+$, it is easy to see that the two states
$(\rho_\pm,u_\pm+\beta t)$ can be connected by a contact
discontinuity $ x=u_{\pm}t+\frac{1}{2}\beta t^{2}$. So the solution
can be expressed as
\begin{equation}\label{2.2}
(\rho,u)(x,t)=\left\{\begin{array}{ll}
(\rho_-,u_-+\beta t),\ \ \ \ \ \ \ \ \ \ \ x<u_{-}t+\frac{1}{2}\beta t^{2},\\
(\rho_+,u_++\beta t),\ \ \ \ \ \ \ \ \ \ \ x>u_{+}t+\frac{1}{2}\beta
t^{2},
\end{array}\right.
\end{equation}

For the case $u_->u_+$,  the characteristics originating from the
origin overlap in the domain  $\{(x,t): u_+t+\frac{1}{2}\beta
t^{2}<x<u_-t+\frac{1}{2}\beta t^{2}\}$, which means that there
exists singularity. A solution containing a weighted
$\delta$-measure supported on a curve will be constructed.

In order to define the measure solution as above, like as in
\cite{Chen-Liu1,Sheng-Zhang}, the two-dimensional weighted
$\delta$-measure $w(t)\delta_S$ supported on a smooth curve
$S=\{(x(s),t(s)):a\leq s\leq b\}$ should be introduced as follows:
\begin{equation}\label{2.3}
\langle
w(s)\delta_S,\psi(s,s)\rangle=\int_a^bw(s)\psi(x(s),t(s))\sqrt{{x'(s)}^2+{t'(s)}^2}ds,
\end{equation}
for any $\psi\in C_0^\infty(R\times R_{+})$.

Let $x=x(t)$ be a discontinuity curve, we consider a piecewise
smooth solution of $(\ref{1.1})$ in the form
\begin{equation}\label{2.4}
(\rho,u)(x,t)=\left\{\begin{array}{ll}
(\rho_-,u_-+\beta t),\ \ \ \ \ \ \ \ \ \ \ \ \ \ \ \ \ \ \ \ \ \ \ x<x(t),\\
(w(t)\delta(x-x(t)),u_\delta(t)),\ \ \ \ \ \ \ \ \ \ x=x(t),\\
(\rho_+,u_++\beta t),\ \ \ \ \ \ \ \ \ \ \ \ \ \ \ \ \ \ \ \ \ \ \
x>x(t),
\end{array}\right.
\end{equation}
in which $u_\delta (t)$ is the assignment of $u$ on this delta shock
wave curve and $u_\delta (t)-\beta t$ is assumed to be a constant.
The delta shock wave solution of the Riemann problem (\ref{1.1}) and
(\ref{1.5}) must obey the following generalized Ranking-Hugoniot
conditions:
\begin{equation}\label{2.5}
\left\{\begin{array}{ll}
\DF{dx(t)}{dt}=\sigma(t)=u_\delta(t),\\
\DF{dw(t)}{dt}= \sigma(t)[\rho]-[\rho u],\\
\DF{d(w(t) u_ \delta(t))}{dt}= \sigma(t)[\rho u]-[\rho u^2]+\beta
w(t),
\end{array}\right.
\end{equation}
and the over-compressive entropy condition
\begin{equation}\label{2.6}
\lambda(\rho_+,u_+)<\sigma(t)<\lambda(\rho_-,u_-),\ \ namely
 \  \ u_++\beta t<u_\delta(t)<u_-+\beta t.
\end{equation}

In  $(\ref{2.5})$, it should be remarkable that $$[\rho
u]=\rho_+(u_++\beta t)-\rho_-(u_-+\beta t), \ \ [\rho
u^{2}]=\rho_+(u_++\beta t)^{2}-\rho_-(u_-+\beta t)^{2}.$$

Through solving $(\ref{2.5})$ with $x(0)=0,\ w(0)=0$, we obtain
\begin{equation}\label{2.7}
\left\{\begin{array}{ll}
u_\delta(t)=\sigma(t)=\sigma_{0}+\beta t,\\
x(t)=\sigma_{0}t+\frac{1}{2}\beta t^{2},\\
w(t)=-\sqrt{\rho_-\rho_+}(u_+-u_-)t,
\end{array}\right.
\end{equation}
with  $\sigma_{0}=
\DF{\sqrt{\rho_-}u_-+\sqrt{\rho_+}u_+}{\sqrt{\rho_-}+\sqrt{\rho_+}}$.

It is easy to prove that the delta shock wave solution $(\ref{2.4})$
with  $(\ref{2.7})$ satisfy the system  $(\ref{1.1})$ in the
distributional sense. That is to say, the following identities
\begin{equation}\label{2.8}
\left\{\begin{array}{ll} \langle\rho,\psi_t\rangle+\langle\rho
u,\psi_x\rangle=0,
\\
\langle\rho u,\psi_t\rangle+\langle\rho
u^2,\psi_x\rangle=-\langle\beta\rho,\psi\rangle,
\end{array}\right.
\end{equation}
holds for any test function $\psi\in C_0^\infty(R\times R_{+})$,
 in which
$$\langle\rho u,\psi\rangle=\int_0^\infty\int_{-\infty}^\infty\hat{\rho}_0 \hat{u}_0\psi dxdt+\langle w(t)u_\delta(t)\delta_S,\psi\rangle,$$
with
$$\hat{\rho}_0=\rho_-+[\rho]H(x-\sigma t),\ \ \hat{u}_{0}=u_{-}-\beta t+[u]H(x-\sigma t).$$

From the above discussions, we can concluded that the Riemann
problem (\ref{1.1}) and (\ref{1.5}) can be solved by three kinds of
solutions: one contact discontinuity, two contact discontinuities
with the vacuum state between them (see Fig.2.1), or the delta shock
wave (see Fig.2.2) connecting two states $(\rho_{\pm},u_{\pm}+\beta
t)$.

%%%%%%%%%%%%%%%%%%%%%%%%%%%%%%%% Fig.2.1%%%%%%%%%%%%%%%%%%%%%%%%%
%TeXCAD Options
%\grade{\on}
%\emlines{\off}
%\beziermacro{\on}
%\reduce{\on}
%\snapping{\off}
%\quality{8}
%\graddiff{.01}
%\snapasp{1}
%\zoom{4.0000}
\unitlength 1mm % = 2.85pt
\linethickness{0.4pt}
\ifx\plotpoint\undefined\newsavebox{\plotpoint}\fi % GNUPLOT compatibility
\begin{picture}(169,69.25)(10,0)
%\vector(15,18)(82,18)
\put(82,18){\vector(1,0){.07}} \put(15,18){\line(1,0){67}}
%\end
%\vector(9,20.25)(9,68)
\put(9,68){\vector(0,1){.07}} \put(9,20.25){\line(0,1){47.75}}
%\end
\qbezier(40,18)(49.75,47.63)(75.5,57.25)
%\qbezier(40,18)(16.13,34.75)(22.5,65.75)
%\qbezier(40,18)(15.13,34.75)(18.5,65.75)
\qbezier(40,18)(17.13,34.75)(26.5,65.75)
%\qbezier(40,18)(19.13,34.75)(31.5,65.75)
%\qbezier(40,18)(21.13,34.75)(34.5,65.75)
%\vector(96,18.5)(166.25,18.25)
\put(166.25,18.25){\vector(1,0){.07}}
\put(96,18.5){\line(1,0){68.78125}}
%\end
%\vector(88.25,21)(88.25,68)
\put(91.25,68){\vector(0,1){.07}} \put(91.25,21){\line(0,1){47}}
%\end
\qbezier(125,18.5)(126.63,54.25)(144.75,68)
\qbezier(125,18.5)(138.75,44.63)(161.25,55.5) \put(5.5,65.25){$t$}
\put(84,18){$x$} \put(88,65){$t$} \put(169,18){$x$}
\put(19,64.25){$J_1$} \put(77,59.5){$J_2$} \put(39.75,14.75){0}
\put(125,15){0} \put(136.25,66.75){$J_1$} \put(164.5,57){$J_2$}
\put(10,23){$(\rho_-,u_-+\beta t)$}
\put(53.25,32.75){$(\rho_+,u_++\beta t)$} \put(32,46){Vac.}
\put(100,40){$(\rho_-,u_-+\beta t)$} \put(134,51.75){Vac.}
\put(145,32){$(\rho_+,u_++\beta t)$}
\put(45,11.5){\makebox(0,0)[cc]{(a)$\ \ u_-<0<u_+$}}
\put(130,12.75){\makebox(0,0)[cc]{(b) $\ \ 0<u_-<u_+$}}
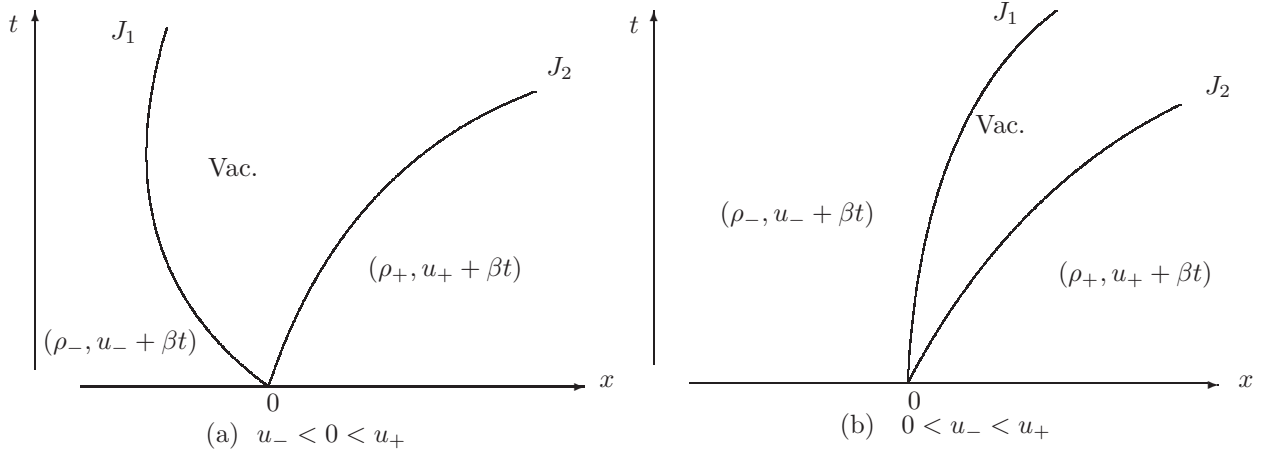
\put(90,6){\makebox(0,0)[cc] {Fig.2.1 The Riemann solution to (1.1)
and (1.5) when $\beta>0$.}}

\end{picture}

%%%%%%%%%%%%%%%%%%%%%%%%%%%%%%% Fig.2.2%%%%%%%%%%%%%%%%%%%%%%%%%%%%%%
%TeXCAD Options
%\grade{\on}
%\emlines{\off}
%\beziermacro{\on}
%\reduce{\on}
%\snapping{\off}
%\quality{8}
%\graddiff{.01}
%\snapasp{1}
%\zoom{4.0000}
\unitlength 1mm % = 2.85pt
\linethickness{0.5pt}
\ifx\plotpoint\undefined\newsavebox{\plotpoint}\fi % GNUPLOT compatibility
\begin{picture}(182,63)(10,0)
%\vector(12,12.5)(85.5,12.5)
\put(85.5,12.5){\vector(1,0){.07}} \put(12,12.5){\line(1,0){73.5}}
%\end
%\vector(4.75,15)(4.75,61)
\put(4.75,61){\vector(0,1){.07}} \put(4.75,15){\line(0,1){46}}
%\end
\qbezier(38.25,12.5)(39.75,40.25)(67.25,55.25)
%\vector(104,13.5)(176,13.5)
\put(176,13.5){\vector(1,0){.07}} \put(104,13.5){\line(1,0){72}}
%\end
%\vector(94.75,14.5)(94.75,63)
\put(94.75,63){\vector(0,1){.07}} \put(94.75,14.5){\line(0,1){48.5}}
%\end
\put(1.75,58.75){$t$} \put(87.75,11.75){$x$} \put(66,56){$\delta
\!S$} \put(92,61.75){$t$} \put(177,13){$x$} \put(110,50){$\delta
\!S$} \put(21,35.75){$(\rho_-,u_-+\beta t)$}
\put(55,35.75){$(\rho_+,u_++\beta t)$}
\put(106.75,35){$(\rho_-,u_-+\beta t)$}
\put(150,35.25){$(\rho_+,u_++\beta t)$} \put(38,8.5){0}
\put(115,9){0} \put(45,8.5){(a)$\ \ \beta>0$} \put(127.5,8.5){(b)$\
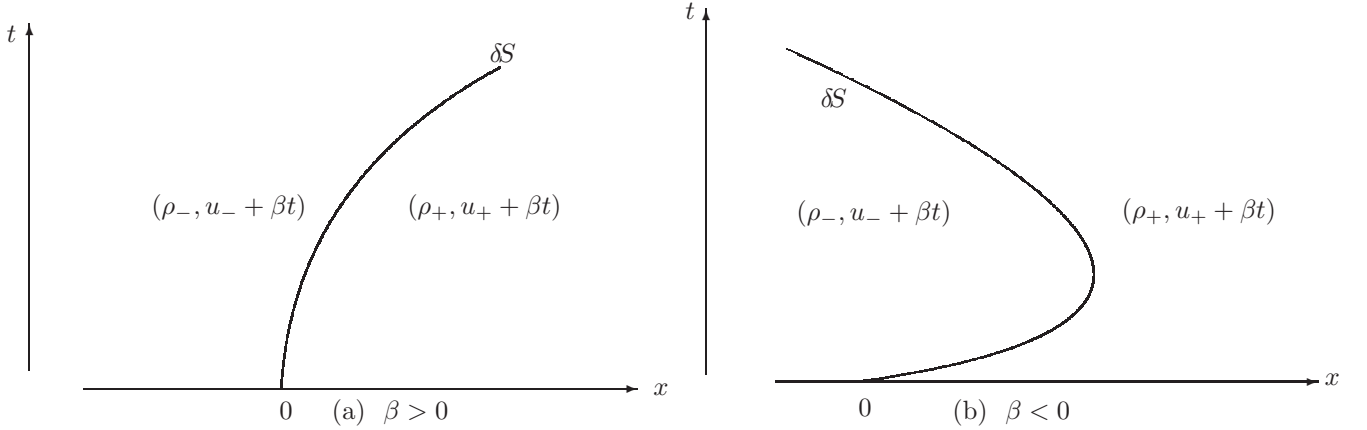
\ \beta<0$} \put(90,5){\makebox(0,0)[cc] {Fig.2.2 The delta shock
wave solution to (1.1) and (1.5) when $u_+<u_-$ and $\sigma_0>0$. }}

 \qbezier(115,13.5)(182,22.75)(105.5,57.75)
\end{picture}

\section{Riemann problem for a modified conservative system of (\ref{1.2})}
In this section, we are devoted to the study of the Riemann problem
for a conservative system of (\ref{1.2}) in detail. Let us introduce
the new velocity $v(x,t)=u(x,t)-\beta t$, then the system
(\ref{1.2}) can be reformulated into a conservative form as follows:

\begin{equation}\label{3.1}
\left\{\begin{array}{ll}
\rho_t+(\rho (v+\beta t))_x=0,\\
(\rho( v+P))_t+(\rho (v +P)(v+\beta t))_x=0.
 \end{array}\right.
\end{equation}
In fact, the change of variable was introduced by Faccanoni and
Mangeney \cite{Faccanoni-Mangeney} to study the shock and
rarefaction waves of the Riemann problem for the shallow water
equations with a with Coulomb-like friction. Here, we use this
transformation to study the delta shock wave for the system
(\ref{1.2}).

Now we want to deal with the Riemann problem for the conservative
system (\ref{3.1}) with the same
 Riemann initial data (\ref{1.5}) as follows:
\begin{equation}\label{3.2}
(\rho,v)(x,0)=\left\{\begin{array}{ll}
(\rho_-,u_-),\ \ x<0,\\
(\rho_ +,u_+),\ \ x>0.
\end{array}\right.
\end{equation}
We shall see hereafter that the Riemann solutions to (\ref{1.2}) and
(\ref{1.5}) can be obtained immediately from the Riemann solutions
to (\ref{3.1}) and (\ref{3.2}) by using the transformation of state
variables $(\rho,u)(x,t)=(\rho,v+\beta t)(x,t)$.

The system (\ref{3.1}) can be rewritten in the quasi-linear form
\begin{equation}\label{3.3}
\left(\begin{array}{lll}
1 & 0 \\
v+P+\rho P' & \rho
\end{array}\right)
\left(\begin{array}{lll}
\rho  \\
v
\end{array}\right)_{t}
+\left(\begin{array}{lll}
v+\beta t & \rho  \\
(v+P+\rho P')(v+\beta t) &\rho(2v+\beta t+P)
\end{array}\right)
\left(\begin{array}{lll}
\rho  \\
v
\end{array}\right)_{x}
=\left(\begin{array}{lll}
0  \\
0
\end{array}\right).
\end{equation}
It can be derived directly from (\ref{3.3}) that the conservative
system (\ref{3.1}) has two eigenvalues
$$\lambda_1(\rho,v)=v+\beta t-\frac{A\alpha}{\rho^{\alpha}},\ \ \lambda_2(\rho,v)=v+\beta t,$$
whose corresponding right eigenvectors can be expressed respectively
by
$$r_1=(\rho,-\frac{A\alpha}{\rho^{\alpha}})^T,\ \ r_2=(1,0)^T.$$
So (\ref{3.1}) is strictly hyperbolic for $\rho>0$. Moreover,
$\bigtriangledown\lambda_1\cdot r_1\neq0$ and
$\bigtriangledown\lambda_2\cdot r_2=0$. Then it can be concluded
that $\lambda_1$ is genuinely nonlinear whose associated waves are
shock waves denoted by $S_1$ or rarefaction waves denoted by $R_1$,
see \cite{Smoller}. Then the Riemann invariants along the
characteristic fields may be chosen as
$$w=v-\frac{A}{\rho^{\alpha}},\ \ z=v,$$ which should satisfy $\bigtriangledown w\cdot r_1=0$ and
$\bigtriangledown z\cdot r_2=0$, respectively.

    Let us draw our attention on the elementary waves for the system
(\ref{3.1}) in detail. We first consider the rarefaction wave which
is a one-parameter family of states connecting a given state. This
kind of continuous solution satisfying the system (\ref{3.1}) can be
obtained by determining the integral curves of the first
characteristic fields. It is worthwhile to notice that the 1-Riemann
invariant is conserved in the 1-rarefaction wave.

For a given left state $(\rho_{-},u_{-})$, the 1-rarefaction wave
curve $R_{1}(\rho_{-},v_{-})$ in the phase plane which is the set of
states connected on the right, should satisfy

\begin{equation}\label{3.4}
R_{1}(\rho_{-},u_{-}): \left\{\begin{array}{ll}
\frac{dx}{dt}=\lambda_{1}(\rho,v)=v+\beta t-\frac{A\alpha}{\rho^{\alpha}},\\
v-\frac{A}{\rho^{\alpha}}=u_{-}-\frac{A}{\rho_{-}^{\alpha}}=w_{-},\\
\lambda_{1}(\rho_{-},u_{-})\leq \lambda_{1}(\rho,v).
\end{array}\right.
\end{equation}

  By differentiating $v$ with respect to $\rho$ in the second
equation in (\ref{3.4}), we have
$$\frac{dv}{d\rho}=-\frac{A\alpha}{\rho^{\alpha+1}}<0,$$
$$\frac{d^{2}v}{d\rho^{2}}=\frac{A\alpha(\alpha+1)}{\rho^{\alpha+2}}>0.$$
Thus, the 1-rarefaction wave is made up of the half-branch of
$R_{1}(\rho_{-},u_{-})$ satisfying $v\geq u_{-}$ and $\rho\leq
\rho_{-}$, which is convex in the $(\rho,v)$ plane.

  Let us compute the solution $(\rho,v)$ at a point in the interior
of the 1-rarefaction wave, then it follows from the first equation
in (\ref{3.4}), we have
\begin{equation}\label{3.5}
v-\frac{A\alpha}{\rho^{\alpha}}=\frac{x}{t}-\beta t.
\end{equation}
By combining (\ref{3.5}) with the second equation in (\ref{3.4}), we
get
\begin{equation}\label{3.6}
(\rho,v)(x,t)=\big(\big(\frac{A(1-\alpha)}{\frac{x}{t}-\beta
t-w_{-}}\big)^{\frac{1}{\alpha}},\frac{\frac{x}{t}-\beta t-\alpha
w_{-}}{1-\alpha}\big).
\end{equation}

Let us return our attention on the shock wave which is a piecewise
constant discontinuous solution, satisfying the Rankine-Hugoniot
conditions and the entropy condition. Here  the Ranking-Hugoniot
conditions can be derived in a standard method as in \cite{Smoller},
since the parameter $t$ only appears
 in the flux functions in the conservative
system (\ref{3.1}). For a bounded discontinuity at $x=x(t)$, let us
denote $\sigma(t)=x'(t)$, then the Ranking-Hugoniot conditions for
the conservative system (\ref{3.1}) can be expressed as
\begin{equation}\label{3.7}
\left\{\begin{array}{ll}
-\sigma(t)\rho+[\rho (v+\beta t)]=0,\\
-\sigma(t)[\rho (v+P)]+[\rho(v+P) (v+\beta t)]=0,
\end{array}\right.
\end{equation}
where $[\rho]=\rho_{r}-\rho_l$ with
$\rho_{l}=\rho(x(t)-0,t)$,\ $\rho_{r}=\rho(x(t)+0,t)$, in which
$[\rho]$ denote the jump of $\rho$ across the discontinuity, etc. It
is clear that the propagation speed of the discontinuity depends on
the parameter $t$, which is obviously different from classical
hyperbolic conservation laws.

If $\sigma(t)\neq0$, then it follows from (\ref{3.7}) that

\begin{equation}\label{3.8}
\rho_{r}\rho_{l}(v_{r}-v_{l})((v_{r}-\frac{A}{\rho_{r}^{\alpha}})-(v_{l}-\frac{A}{\rho_{l}^{\alpha}}))=0,
\end{equation}
from which we have $v_{r}=v_{l}$ or
$v_{r}-\frac{A}{\rho_{r}^{\alpha}}=v_{l}-\frac{A}{\rho_{l}^{\alpha}}$.

Thus, for a given left state $(\rho_{-},u_{-})$, with the latex
entropy condition in mind, the 1-shock wave curve
$S_{1}(\rho_{-},u_{-})$ in the $(\rho,v)$ plane which is the set of
states connected on the right, should satisfy
\begin{equation}\label{3.9}
S_{1}(\rho_{-},u_{-}): \left\{\begin{array}{ll}
\sigma_{1}(t)=\frac{\rho v-\rho u_{-}}{\rho-\rho_{-}}+\beta t,\\
v-\frac{A}{\rho^{\alpha}}=u_{-}-\frac{A}{\rho_{-}^{\alpha}}=w_{-},\\
\rho>\rho_{-}, \ \ v<u_{-},
\end{array}\right.
\end{equation}
which indicates the 1-rarefaction wave and 1-shock wave are
different branch of the same curve.

Moreover, from (\ref{3.8}),  for a given left state
$(\rho_{-},u_{-})$, the 2-contact discontinuity curve
$J(\rho_{-},u_{-})$ in the $(\rho,v)$ plane which is the set of
states connected on the right, should satisfy

\begin{equation}\label{3.10}
J(\rho_{-},u_{-}):\ \  \ \sigma(t)=v+\beta t=u_{-}+\beta t.
\end{equation}

%%%%%%%%%%%%%%%%%%%%%%%%%%%%%%%% FIG.1%%%%%%%%%%%%%

%TeXCAD Options
%\grade{\on}
%\emlines{\off}
%\beziermacro{\on}
%\reduce{\on}
%\snapping{\off}
%\quality{8}
%\graddiff{.01}
%\snapasp{1}
%\zoom{4.0000}
\unitlength 0.9mm % = 2.85pt
\linethickness{0.4pt}
\ifx\plotpoint\undefined\newsavebox{\plotpoint}\fi % GNUPLOT compatibility
\begin{picture}(200,69)(-15,0)
%\vector(15.75,13.75)(95.75,13.75)
\put(95.75,13.75){\vector(1,0){.07}}
\put(15.75,13.75){\line(1,0){80}}
%\end
%\vector(15.75,14)(15.75,69)
\put(5.75,69){\vector(0,1){.07}} \put(5.75,14){\line(0,1){55}}
%\end
%\emline(35.25,67.75)(35.25,13.75)
%\put(35.25,67.75){\line(0,-1){54}}
\bezier{50}(35.25,67.75)(35.25,23.25)(35.25,13.75)
%\end
%\emline(65,68)(65,14.25)
\put(65,68){\line(0,-1){54}}
%\end
\qbezier(37.5,65.5)(46,23.25)(88.5,21) \put(12.5,67.5){$\rho$}
\put(99,13.75){$v$} %\put(15.5,9.75){0}
\put(26,41.25){I\!I\!I} \put(55,46.25){I\!I} \put(76.5,45.5){I}
\put(44.75,53.75){$S_1$} \put(69,55.25){$J$} \put(80,23.25){$R_1$}
\put(31,50.75){$S_{\delta}$}
\put(30,9.5){$u_--\frac{A}{\rho_-^\alpha}$} \put(63.5,11){$u_-$}
\put(67,30){($\rho_-$,$u_-$)} \put(60.75,5.5){\makebox(0,0)[cc]
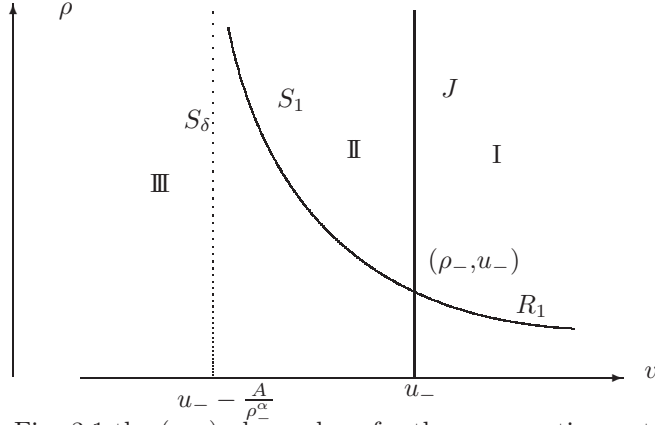
{Fig. 3.1  the $(\rho,v)$ phase plane for the conservative system
(\ref{3.1}).}}
\end{picture}

Let us now consider the Riemann problem (\ref{3.1}) and (\ref{3.2}).
In the $(\rho,v)$ phase plane, for a given left state
$(\rho_{-},u_{-})$, the set of  states connected on the right
consist of the 1-rarefaction wave $R_{1}(\rho_{-},u_{-})$, the
1-shock wave $S_{1}(\rho_{-},u_{-})$ and the 2-contact discontinuity
curve $J(\rho_{-},u_{-})$. It is clear to see that
$R_{1}(\rho_{-},u_{-})$ has the line
$S_{\delta}:v=u_{-}-\frac{A}{\rho_{-}^{\alpha}}$ and
$S_{1}(\rho_{-},u_{-})$ has the positive $v$-axis as their
asymptotic lines, respectively.

In view of the right state $(\rho_{+},u_{+})$ in different
positions, one wants to construct the unique global Riemann solution
of (\ref{3.1}) and (\ref{3.2}). However, as in \cite{Guo}, if
$u_{+}\leq u_{-}-\frac{A}{\rho_{-}^{\alpha}}$ is satisfied, the
Riemann solution of (\ref{3.1}) and (\ref{3.2}) can not be
constructed by using only the elementary waves including shocks,
rarefaction waves and contact discontinuities. In this nonclassical
situation, the concept of delta shock wave should be introduced such
as in \cite{Guo,Guo-Sheng-Zhang,Wang} and be discussed later.

Draw all the curves  $R_{1}(\rho_{-},u_{-})$,
$S_{1}(\rho_{-},u_{-})$
 $J(\rho_{-},u_{-})$ and $S_{\delta}$ in the the $(\rho,v)$ phase
 plane, thus the phase plane is divided into three regions  I,
I\!I and I\!I\!I (See Fig.3.1), where
$${\rm I}=\{(\rho,v)|v\geq u_-\},$$
$${\rm I\!I}=\{(\rho,v)|u_{-}-\frac{A}{\rho_{-}^\alpha}<v< u_-\},$$
$${\rm I\!I\!I}=\{(\rho,v)|v\leq u_--\frac{A}{\rho_{-}^\alpha}\}.$$
According to the right state  $(\rho_{+},u_{+})$ in different
regions, the unique  global Riemann solution of (\ref{3.1}) and
(\ref{3.2}) can be constructed connecting two constant states
$(\rho_{-},u_{-})$ and $(\rho_{+},u_{+})$

 If
$(\rho_+,u_+)\in$ I, namely $u_+>u_{-}$, then the Riemann solution
consists of 1-rarefaction wave $R_1$ and a 2-contact discontinuity
$J$ with an intermediate constant state $(\rho_\ast,v_\ast)$
determined uniquely
 by
\begin{equation}\label{3.11}
\left\{\begin{array}{ll}
v_{\ast}-\frac{A}{\rho_{\ast}^{\alpha}}=u_--\frac{A}{\rho_-^{\alpha}}=w_{-},\\
u_+=v_{\ast}.
\end{array}\right.
\end{equation}
which immediately leads to

\begin{equation}\label{3.12}
(\frac{A}{\rho_\ast^{\alpha}},v_\ast)=(u_+-u_-+\frac{A}{\rho_-^{\alpha}},u_+),
\end{equation}
or
\begin{equation}\label{3.13}
(\rho_\ast,v_\ast)=\big(\big(\frac{A}{(u_+-u_-+\frac{A}{\rho_-^{\alpha}}}\big)^\frac{1}{\alpha},u_+\big),
\end{equation}
Thus, the Riemann solution of (\ref{3.1}) and (\ref{3.2}) can be
express as
\begin{equation}\label{3.14}
(\rho,v)(x,t)=\left\{\begin{array}{ll}
(\rho_{-},u_{-}+\beta t),\ \ x<x_{1}^{-}(t),\\
R_{1},\ \ \ \ \ \ \ \ \ \ \ \ \ \ \ x_{1}^{-}(t)<x<x_{1}^{+}(t),\\
(\rho_{\ast},u_{\ast}+\beta t),\ \ \ x_{1}^{+}(t)<x<x_{2}(t),\\
(\rho_{+},u_{+}+\beta t),\ \ x>x_{2}(t),\\
\end{array}\right.
\end{equation}
in which
\begin{equation}\label{3.15}
x_{1}^{-}(t)=(u_--\frac{A}{\rho_-^{\alpha}})t+\frac{1}{2}\beta
t^{2},\ \
x_{1}^{+}(t)=(u_\ast-\frac{A}{\rho_\ast^{\alpha}})t+\frac{1}{2}\beta
t^{2},\end{equation}
\begin{equation}\label{3.16}
x_{2}(t)=u_{+}t+\frac{1}{2}\beta t^{2},\end{equation} and the state
$(\rho_1,u_1)$ in $R_{1}$ can be calculated by (\ref{3.6}).

If $(\rho_+,u_+)\in$ I\!I, namely
$u_{-}-\frac{A}{\rho_{-}^\alpha}<u_+< u_-$, then the Riemann
solution consists of 1-shock wave $S_1$ and a 2-contact
discontinuity $J$ with an intermediate constant state
$(\rho_\ast,v_\ast)$ determined uniquely by (\ref{3.13}). Thus, the
Riemann solution of (\ref{3.1}) and (\ref{3.2}) can be express as
\begin{equation}\label{3.17}
(\rho,v)(x,t)=\left\{\begin{array}{ll}
(\rho_{-},u_{-}+\beta t),\ \ x<x_{1}(t),\\
(\rho_{\ast},u_{\ast}+\beta t),\ \ \ x_{1}(t)<x<x_{2}(t),\\
(\rho_{+},u_{+}+\beta t),\ \ x>x_{2}(t),\\
\end{array}\right.
\end{equation}
in which the position of $S_1$ is given by
\begin{equation}\label{3.18}
x_{1}(t)=\frac{\rho_{\ast}v_{\ast}-\rho_-u_-}{\rho_{\ast}-\rho_-}t+\frac{1}{2}\beta
t^{2},
\end{equation}
 and $x_{2}(t)$ is given by (\ref{3.16}).

On the other hand, when $(\rho_+,u_+)\in$ I\!I\!I, namely $
u_{+}\leq u_{-}-\frac{A}{\rho_{-}^{\alpha}}$, then there exist a
nonclassical situation where the Cauchy problem does not own a weak
$L^\infty$-solution. In order to solve the Riemann problem
(\ref{3.1}) and (\ref{3.2}) in the framework of nonclassical
solution, a solution containing a weighted $\delta$-measure
supported on a curve should be defined such as in
\cite{Chen-Liu1,Pan-Han,Sheng-Zhang}. In what follows, let us
provide the definition of delta shock wave solution to the Riemann
problem (\ref{3.1}) and (\ref{3.2}). Let us also refer to
\cite{Danilvo-Shelkovich1,Danilvo-Shelkovich2,Kalisch-Mitrovic1,Kalisch-Mitrovic2}
about the more exact definition of generalized delta shock wave
solution for related systems with delta measure initial data.

\begin{defn}\label{defn:3.1}
Let $(\rho,v)$ be a pair of distributions in which $\rho$ has the
form of
\begin{equation}\label{3.19}
\rho(x,t)=\hat{\rho}(x,t)+w(x,t)\delta_{S},
\end{equation}
in which $\hat{\rho},v\in L^{\infty}(R\times R_{+})$. Then,
$(\rho,v)$ is called as the delta shock wave solution to the Riemann
problem (\ref{3.1}) and (\ref{3.2}) if it satisfies
\begin{equation}
\left\{\begin{array}{ll}\nonumber
\langle\rho,\psi_{t}\rangle+\langle\rho (v+\beta t),\psi_{x}\rangle=0,\\
\langle\rho( v+P)),\psi_{t}\rangle+\langle\rho ( v+P)(v+\beta
t)),\psi_{x}\rangle=0,
 \end{array}\right.
\end{equation}
for any $\psi\in C_0^\infty(R\times R^{+})$. Here we take
 $$\langle\rho
(v+P)(v+\beta
t)),\psi\rangle=\int_{0}^{\infty}\int_{-\infty}^{\infty}(\widehat{\rho}(v-\frac{A}{\widehat{\rho}^\alpha})(v+\beta
t))\psi dxdt+\langle w(t)v_{\delta}(t)(v_{\delta}(t)+\beta
t)\delta_{S},\psi\rangle,$$ as an example to explain the inner
product, in which we use the symbol $S$ to express the smooth curve
with the Dirac delta function supported on it, $v_\delta$ is the
value of $v$ and $\frac{A}{\rho^\alpha}$ is equal to zero on this
delta shock wave $S$.
\end{defn}

With the above definition, if $(\rho_+,u_+)\in$ I\!I\!I and
$u_{+}<u_{-}-\frac{A}{\rho_{-}^\alpha}$, a piecewise smooth solution
of the Riemann problem (\ref{3.1}) and (\ref{3.2}) should be
introduced in the form
\begin{equation}\label{3.20}
(\rho,v)(x,t)=\left\{\begin{array}{ll}
(\rho_-,u_-),\ \ \ \ \ \ \ \ \ \ \ \ \ \ \ \ \ \ \ \ x<x(t),\\
(w(t)\delta(x-x(t)),v_{\delta}),\ \ \  \ x=x(t),\\
(\rho_+,u_+),\ \ \ \ \ \ \ \ \ \ \ \ \ \ \ \ \ \ \ \ x>x(t),
\end{array}\right.
\end{equation}
where $x(t)$, $w(t)$ and $\sigma(t)=x'(t)$ denote respectively the
location, weight and propagation speed of the delta shock, and
$v_{\delta}$ indicates the assignment of $v$ on this  delta shock
wave. It is remarkable that the value of $v$ should be given on the
delta shock curve $x=x(t)$ such that the product of $\rho$ and $v$
can be defined in the sense of distributions. When
$u_{+}=u_{-}-\frac{A}{\rho_{-}^\alpha}$, it can be discussed
similarly and we omit it.

The delta shock wave solution of the form (\ref{3.20}) to the the
Riemann problem (\ref{3.1}) and (\ref{3.2}) should obey the
generalized Rankine-Hugoniot conditions
\begin{equation}\label{3.21}
\left\{\begin{array}{ll}
\DF{dx(t)}{dt}=\sigma(t)=v_{\delta}+\beta t,\\[4pt]
\DF{dw(t)}{dt}=\sigma(t)[\rho]-[\rho (v+\beta t)],\\[4pt]
\DF{d(w(t)v_{\delta})}{dt}=\sigma(t)[\rho(
v-\frac{A}{\rho^\alpha})]-[\rho( v-\frac{A}{\rho^\alpha}) (v+\beta
t)],
\end{array}\right.
\end{equation}
with initial data $x(0)=0$ and $w(0)=0$. In addition, for the unique
solvability of the above Cauchy problem, it is necessary to require
that the value of $v_{\delta}$ to be a constant along the trajectory
of delta shock wave (see \cite{Danilvo-Shelkovich2} for details).
The derivation process of the generalized Rankine-Hugoniot
conditions is similar to that in \cite{Shen1,Shen2,Sun} and we omit
it here. In order to ensure the uniqueness of Riemann solutions, an
over-compressive entropy condition for the delta shock wave should
be proposed by
\begin{equation}\label{3.22}
\lambda_{1}(\rho_{+},u_{+})<\lambda_{2}(\rho_{+},u_{+})<\sigma(t)<\lambda_{1}(\rho_{-},u_{-})<\lambda_{2}(\rho_{-},u_{-}),
\end{equation}
such that we have
\begin{equation}\label{3.23}
u_{+}<v_{\delta}<u_{-}-\frac{A}{\rho_{-}^\alpha},
\end{equation}
which implies that all the characteristics on both sides of the
delta shock are in-coming.

It follows from (\ref{3.21}) that
\begin{equation}\label{3.24}
\DF{dw(t)}{dt}=v_{\delta}(\rho_{+}-\rho_{-})-(\rho_{+}u_{+}-\rho_{-}u_{-}),
\end{equation}
\begin{equation}\label{3.25}
v_{\delta}\DF{dw(t)}{dt}=v_{\delta}\big((\rho_{+}u_{+}-\rho_{-}u_{-})-(\frac{A}{\rho_{+}^{\alpha-1}}-\frac{A}{\rho_{-}^{\alpha-1}})\big)-
(\rho_{+}u_{+}^{2}-\rho_{-}u_{-}^{2})+\big(\frac{Au_{+}}{\rho_{+}^{\alpha-1}}-\frac{Au_{-}}{\rho_{-}^{\alpha-1}}\big),
\end{equation}
Thus, we have
\begin{equation}\label{3.26}
(\rho_{+}-\rho_{-})v_{\delta}^{2}-\big(2(\rho_{+}u_{+}-\rho_{-}u_{-})-(\frac{A}{\rho_{+}^{\alpha-1}}-\frac{A}{\rho_{-}^{\alpha-1}})\big)v_{\delta}+
(\rho_{+}u_{+}^{2}-\rho_{-}u_{-}^{2})-\big(\frac{Au_{+}}{\rho_{+}^{\alpha-1}}-\frac{Au_{-}}{\rho_{-}^{\alpha-1}}\big)=0,
\end{equation}

For convenience, let us denote
\begin{equation}\label{3.27}
w_{0}=\sqrt{\rho_{+}\rho_{-}(u_{+}-u_{-})\big((u_{+}-u_{-})-(\frac{A}{\rho_{+}^\alpha}-\frac{A}{\rho_{-}^\alpha})\big)+
\frac{1}{4}\big(\frac{A}{\rho_{+}^{\alpha-1}}-\frac{A}{\rho_{-}^{\alpha-1}}\big)^2}-
\frac{1}{2}\big(\frac{A}{\rho_{+}^{\alpha-1}}-\frac{A}{\rho_{-}^{\alpha-1}}\big)>0,
\end{equation}

If $\rho_+\neq\rho_-$, with the entropy condition (\ref{3.22}) in
mind, one can obtain directly from (\ref{3.27}) that
\begin{equation}\label{3.28}
v_{\delta}=\frac{\rho_+ u_+-\rho_-u_-+w_{0}}{\rho_+-\rho_-},
\end{equation}
which enables us to get
\begin{equation}\label{3.29}
\sigma(t)=v_{\delta}+\beta t,\ \ x(t)=v_{\delta}t+\frac{1}{2}\beta
t^{2},\ \ w(t)=w_{0}t.
\end{equation}

Otherwise, if $\rho_+=\rho_-$, then we have

\begin{equation}\label{3.30}
v_{\delta}=\frac{1}{2}(u_++u_--\frac{A}{\rho_{-}^\alpha}).
\end{equation}
In this particular case, we can also get
\begin{equation}\label{3.31}
\sigma(t)=\frac{1}{2}(u_++u_--\frac{A}{\rho_{-}^\alpha})+\beta t,\
x(t)=\frac{1}{2}(u_++u_--\frac{A}{\rho_{-}^\alpha})t+\frac{1}{2}\beta
t^{2},\  w(t)=(\rho_-u_--\rho_+ u_+)t.
\end{equation}

\section{Riemann problem for the approximated system (\ref{1.2})}
In this section, let us return to the Riemann problem (\ref{1.2})
and (\ref{1.5}). If $(\rho_+,u_+)\in$ I, the Riemann solutions to
(\ref{1.2}) and (\ref{1.5})  $R_1+J$ can be represented as
\begin{equation}\label{4.1}
(\rho,u)(x,t)=\left\{\begin{array}{ll}
(\rho_-,u_-+\beta t),\ \ \ \ \ \ \ \ \ \ \ x<x_{1}^-(t),\\
(\rho_1,v_1+\beta t),\ \ \  \ \ \ \ \ \ \ \ \ x_{1}^-(t)<x<x_{1}^+(t),\\
(\rho_*,v_*+\beta t),\ \ \  \ \ \ \ \ \ \ \ \ x_{1}^+(t)<x<x_{2}(t),\\
(\rho_+,u_++\beta t),\ \ \ \ \ \ \ \ \ \ \ x>x_{2}(t),
\end{array}\right.
\end{equation}
where $x_{1}^-(t),\ x_{1}^+(t)$ and $x_{2}(t)$ are given by
(\ref{3.15}) and (\ref{3.16}) respectively, and the states
$(\rho_1,v_1)$ and $(\rho_*,v_*)$ can be calculated as (\ref{3.6})
and (\ref{3.13}). Let us use Fig.4.1(a) to illustrate this situation
in detail, where all the characteristics in the rarefaction wave
fans $R_1$ and contact discontinuity curve $J$ are curved into
parabolic shapes.

If $(\rho_+,u_+)\in$ I\!I, the Riemann solutions to (\ref{1.2}) and
(\ref{1.5}) $S_1+J$ can be represented as
\begin{equation}\label{4.2}
(\rho,u)(x,t)=\left\{\begin{array}{ll}
(\rho_-,u_-+\beta t),\ \ \ \ \ \ \ \ \ \ \ x<x_{1}(t),\\
(\rho_*,v_*+\beta t),\ \ \  \ \ \ \ \ \ \ \ \ x_{1}(t)<x<x_{2}(t),\\
(\rho_+,u_++\beta t),\ \ \ \ \ \ \ \ \ \ \ x>x_{2}(t),
\end{array}\right.
\end{equation}
where $x_{1}(t)$ and $x_{2}(t)$ are given by (\ref{3.18}) and
(\ref{3.16}) respectively and the states $(\rho_*,v_*)$ can be
calculated as (\ref{3.13}). Let us use Fig.4.1(b) to illustrate this
situation in detail, where both the shock wave curve $S_1$ and the
contact discontinuity curve $J$ are curved into parabolic shapes.

Analogously, if $(\rho_+,u_+)\in$ I\!I\!I, namely $u_{+}\leq
u_{-}-\frac{A}{\rho_{-}^\alpha}$, then we can also define the weak
solutions to the Riemann problem (\ref{1.2}) and (\ref{1.5}) in the
sense of distributions below.

\begin{defn}\label{defn:4.1}
Let $(\rho,u)$ be a pair of distributions in which $\rho$ has the
form of (\ref{3.19}), then it is called as the delta shock wave
solution to the Riemann problem (\ref{1.2}) and (\ref{1.5}) if it
satisfies
\begin{equation}\label{4.3}
\left\{\begin{array}{ll}
\langle\rho,\psi_{t}\rangle+\langle\rho u,\psi_{x}\rangle=0,\\
\langle\rho( u+P)),\psi_{t}\rangle+\langle\rho u(
u+P)),\psi_{x}\rangle=-\langle\beta\rho,\psi\rangle,
 \end{array}\right.
\end{equation}
for any $\psi\in C_0^\infty(R\times R^{+})$, in which
 $$\langle\rho
u(
u+P)),\psi\rangle=\int_{0}^{\infty}\int_{-\infty}^{\infty}(\widehat{\rho}u(u-\frac{A}{\widehat{\rho}^\alpha}))\psi
dxdt+\langle w(t)(u_{\delta}(t))^{2}\delta_{S},\psi\rangle,$$ and
$u_{\delta}(t)$ is the assignment of $u$ on this delta shock wave
curve. \end{defn}

%%%%%%%%%%%%%%%%%%%%%%%%%%%%%%%% Fig.3%%%%%%%%%%%%%%%%%%%%%%%%%
%TeXCAD Options
%\grade{\on}
%\emlines{\off}
%\beziermacro{\on}
%\reduce{\on}
%\snapping{\off}
%\quality{8}
%\graddiff{.01}
%\snapasp{1}
%\zoom{4.0000}
\unitlength 1mm % = 2.85pt
\linethickness{0.4pt}
\ifx\plotpoint\undefined\newsavebox{\plotpoint}\fi % GNUPLOT compatibility
\begin{picture}(169,69.25)(10,0)
%\vector(15,18)(82,18)
\put(82,18){\vector(1,0){.07}}
\put(15,18){\line(1,0){67}}
%\end
%\vector(9,20.25)(9,68)
\put(9,68){\vector(0,1){.07}}
\put(9,20.25){\line(0,1){47.75}}
%\end
\qbezier(40,18)(49.75,47.63)(75.5,57.25)
\qbezier(40,18)(16.13,34.75)(22.5,65.75)
\qbezier(40,18)(15.13,34.75)(18.5,65.75)
\qbezier(40,18)(17.13,34.75)(26.5,65.75)
\qbezier(40,18)(19.13,34.75)(31.5,65.75)
\qbezier(40,18)(21.13,34.75)(34.5,65.75)
%\vector(96,18.5)(166.25,18.25)
\put(166.25,18.25){\vector(1,0){.07}}
\put(96,18.5){\line(1,0){68.78125}}
%\end
%\vector(88.25,21)(88.25,68)
\put(91.25,68){\vector(0,1){.07}}
\put(91.25,21){\line(0,1){47}}
%\end
\qbezier(125,18.5)(126.63,54.25)(144.75,68)
\qbezier(125,18.5)(138.75,44.63)(161.25,55.5) \put(5.5,65.25){$t$}
\put(84,18){$x$} \put(85.75,68.75){$t$} \put(169,18){$x$}
\put(23.5,69.25){$R_1$} \put(77,59.5){$J$} \put(39.75,14.75){0}
\put(125,15){0} \put(136.25,66.75){$S_1$} \put(164.5,56){$J$}
\put(10,23){$(\rho_-,u_-+\beta t)$}
\put(53.25,32.75){$(\rho_+,u_++\beta t)$}
\put(32,46){$(\rho_*,v_*+\beta t)$} \put(100,40){$(\rho_-,u_-+\beta
t)$} \put(134,51.75){$(\rho_*,v_*+\beta t)$}
\put(145,32){$(\rho_+,u_++\beta t)$}
\put(45,11.5){\makebox(0,0)[cc]{(a)$\ \
u_--\frac{A}{\rho_-^\alpha}<u_-<u_+$}}
\put(130,12.75){\makebox(0,0)[cc]{(b) $\ \
u_--\frac{A}{\rho_-^\alpha}<u_+<u_-$}} \put(90,6){\makebox(0,0)[cc]
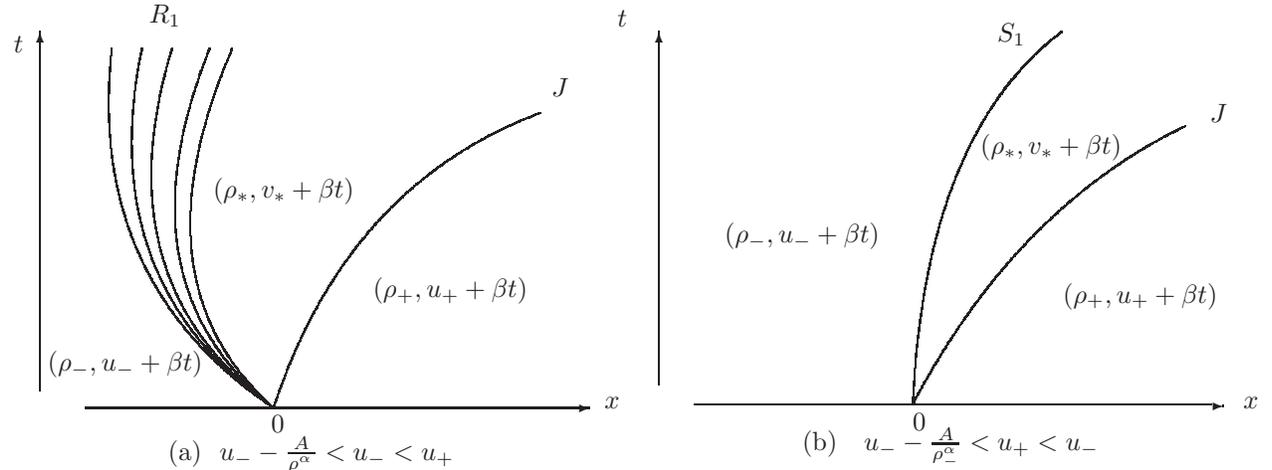
{Fig.4.1 The Riemann solution to (1.2) and (1.5) when
$u_--\frac{A}{\rho_-^\alpha}<u_+$ and $\beta>0$,}}
\put(59,2.5){\makebox(0,0)[cc] {where $(\rho_*,v_*)$ is given by
(3.13).}}

\end{picture}

With the above definition in mind, if $u_{+}<
u_{-}-\frac{A}{\rho_{-}^\alpha}$ is satisfied, then we look for a
piecewise smooth solution to the Riemann problem (\ref{1.2}) and
(\ref{1.5}) in the form
\begin{equation}\label{4.4}
(\rho,u)(x,t)=\left\{\begin{array}{ll}
(\rho_-,u_-+\beta t),\ \ \ \ \ \ \ \ \ \ \ \ \ \ \ \ \ x<x(t),\\
(w(t)\delta(x-x(t)),u_{\delta}(t)),\ \ \  \ x=x(t),\\
(\rho_+,u_++\beta t),\ \ \ \ \ \ \ \ \ \ \ \ \ \ \ \ \ x>x(t),
\end{array}\right.
\end{equation}
It is worthwhile to notice that $u_{\delta}(t)-\beta t$ is assumed
to be a constant based on the result in Sect.2. With the similar
analysis and derivation as before, the delta shock wave solution of
the form (\ref{4.4}) to the Riemann problem (\ref{1.2}) and
(\ref{1.5}) should also satisfy the following generalized
Rankine-Hugoniot conditions
\begin{equation}\label{4.5}
\left\{\begin{array}{ll}
\DF{dx(t)}{dt}=\sigma(t)=u_{\delta}(t),\\[4pt]
\DF{dw(t)}{dt}=\sigma(t)[\rho]-[\rho u],\\[4pt]
\DF{d(w(t)u_{\delta}(t))}{dt}=\sigma(t)[\rho(
u-\frac{A}{\rho^\alpha})]-[\rho u( u-\frac{A}{\rho^\alpha}) ]+\beta
w(t).
\end{array}\right.
\end{equation}
in which the jumps across the discontinuity are
\begin{equation}\label{4.6}
[\rho u]=\rho_+(u_++\beta
t)-\rho_-(u_-+\beta t),
\end{equation}
\begin{equation}\label{4.7}
[\rho u( u-\frac{A}{\rho^\alpha})]=\rho_+(u_++\beta t)(u_++\beta
t-\frac{A}{\rho_{+}^\alpha})-\rho_-(u_-+\beta t)(u_-+\beta
t-\frac{A}{\rho_-^\alpha}).
\end{equation}

In order to ensure the uniqueness to the Riemann problem (\ref{1.2})
and (\ref{1.5}), the over-compressive entropy condition for the
delta shock wave
\begin{equation}\label{4.8}
u_{+}+\beta t<u_{\delta}(t)<u_{-}-\frac{A}{\rho_{-}^\alpha}+\beta t.
\end{equation}
should also be proposed when $u_{+}<
u_{-}-\frac{A}{\rho_{-}^\alpha}$.

Like as before, we can also obtain $x(t),\sigma(t)$ and $w(t)$ from
(\ref{4.5}) and (\ref{4.8}) together. In brief, we have the
following theorem to depict the Riemann solution to (\ref{1.2}) and
(\ref{1.5}) when the Riemann initial data (\ref{1.5}) satisfy $
u_{+}< u_{-}-\frac{A}{\rho_{-}^\alpha}$ and $\rho_{+}\neq\rho_{-}$.

\begin{thm}\label{thm:4.1} If both $ u_{+}< u_{-}-\frac{A}{\rho_{-}^\alpha}$ and
$\rho_{+}\neq\rho_{-}$ are satisfied, then the delta shock solution
to the Riemann solutions to (\ref{1.2}) and (\ref{1.5}) can be
expressed as
\begin{equation}\label{4.9}
\left\{\begin{array}{ll}
\DF{dx(t)}{dt}=\sigma(t)=u_{\delta}(t),\\[4pt]
\DF{dw(t)}{dt}=\sigma(t)[\rho]-[\rho u],\\[4pt]
\DF{d(w(t)u_{\delta}(t))}{dt}=\sigma(t)[\rho(
u-\frac{A}{\rho^\alpha})]-[\rho u( u-\frac{A}{\rho^\alpha}) ]+\beta
w(t).
\end{array}\right.
\end{equation}
in which
\begin{equation}\label{4.10}
\sigma(t)=u_{\delta}(t)=v_{\delta}+\beta t,\ \
x(t)=v_{\delta}t+\frac{1}{2}\beta t^{2},\ \ w(t)=w_{0}t,
\end{equation}
in which $w_{0}$ and $v_{\delta}$ are given by (\ref{3.27}) and
(\ref{3.28}) respectively. \end{thm}

Let us check briefly that the above constructed delta shock wave
solution (\ref{4.9}) and (\ref{4.10}) should satisfy (\ref{1.2}) in
the sense of distributions. The proof of this theorem is completely
analogs to those in \cite {Shen1,Shen2}. Therefore, we only deliver
the main steps for the proof of the second equality in (\ref{4.3})
for completeness. Actually, one can deduce that

\begin{eqnarray}
I&=&\int_{0}^{\infty}\int_{-\infty}^{\infty}(\rho(
u-\frac{A}{\rho^\alpha})\psi_{t}+\rho u(
u-\frac{A}{\rho^\alpha}) \psi_{x})dxdt\nonumber\\
&=&\int_{0}^{\infty}\int_{-\infty}^{x(t)}(\rho_{-}( u_{-}+\beta
t-\frac{A}{\rho_{-}^\alpha})\psi_{t}+\rho_{-}(u_{-}+\beta t)(u
_{-}+\beta
t-\frac{A}{\rho_{-}^\alpha})\psi_{x})dx dt \nonumber\\
&&+\int_{0}^{\infty}\int^{\infty}_{x(t)}(\rho_{+}( u_{+}+\beta
t-\frac{A}{\rho_{+}^\alpha})\psi_{t}+\rho_{+}(u_{+}+\beta t)(u
_{+}+\beta
t-\frac{A}{\rho_{+}^\alpha}) \psi_{x})dx dt \nonumber\\
&&+\int_{0}^{\infty}w_{0}t(v_{\delta}+\beta
t)(\psi_{t}(x(t),t)+(v_{\delta}+\beta
t)\psi_{x}(x(t),t))dt\nonumber.
\end{eqnarray}

It can be derived from (\ref{4.10}) that the curve of delta shock
wave is given by
\begin{equation}\label{4.11}
x(t)=v_{\delta}t+\frac{1}{2}\beta t^{2}.
\end{equation}

%%%%%%%%%%%%%%%%%%%%%%%%%%%%%%% Fig.4%%%%%%%%%%%%%%%%%%%%%%%%%%%%%%
%TeXCAD Options
%\grade{\on}
%\emlines{\off}
%\beziermacro{\on}
%\reduce{\on}
%\snapping{\off}
%\quality{8}
%\graddiff{.01}
%\snapasp{1}
%\zoom{4.0000}
\unitlength 1mm % = 2.85pt
\linethickness{0.5pt}
\ifx\plotpoint\undefined\newsavebox{\plotpoint}\fi % GNUPLOT compatibility
\begin{picture}(182,63)(10,0)
%\vector(12,12.5)(85.5,12.5)
\put(85.5,12.5){\vector(1,0){.07}}
\put(12,12.5){\line(1,0){73.5}}
%\end
%\vector(4.75,15)(4.75,61)
\put(4.75,61){\vector(0,1){.07}}
\put(4.75,15){\line(0,1){46}}
%\end
\qbezier(38.25,12.5)(39.75,40.25)(67.25,55.25)
%\vector(104,13.5)(176,13.5)
\put(176,13.5){\vector(1,0){.07}}
\put(104,13.5){\line(1,0){72}}
%\end
%\vector(94.75,14.5)(94.75,63)
\put(94.75,63){\vector(0,1){.07}}
\put(94.75,14.5){\line(0,1){48.5}}
%\end
\put(1.75,58.75){$t$} \put(87.75,11.75){$x$} \put(66,56){$\delta
\!S$} \put(92,61.75){$t$} \put(177,13){$x$} \put(110,50){$\delta
\!S$} \put(21,35.75){$(\rho_-,u_-+\beta t)$}
\put(55,35.75){$(\rho_+,u_++\beta t)$}
\put(106.75,35){$(\rho_-,u_-+\beta t)$}
\put(150,35.25){$(\rho_+,u_++\beta t)$} \put(38,8.5){0}
\put(115,9){0} \put(45,8.5){(a)$\ \ \beta>0$} \put(127.5,8.5){(b)$\
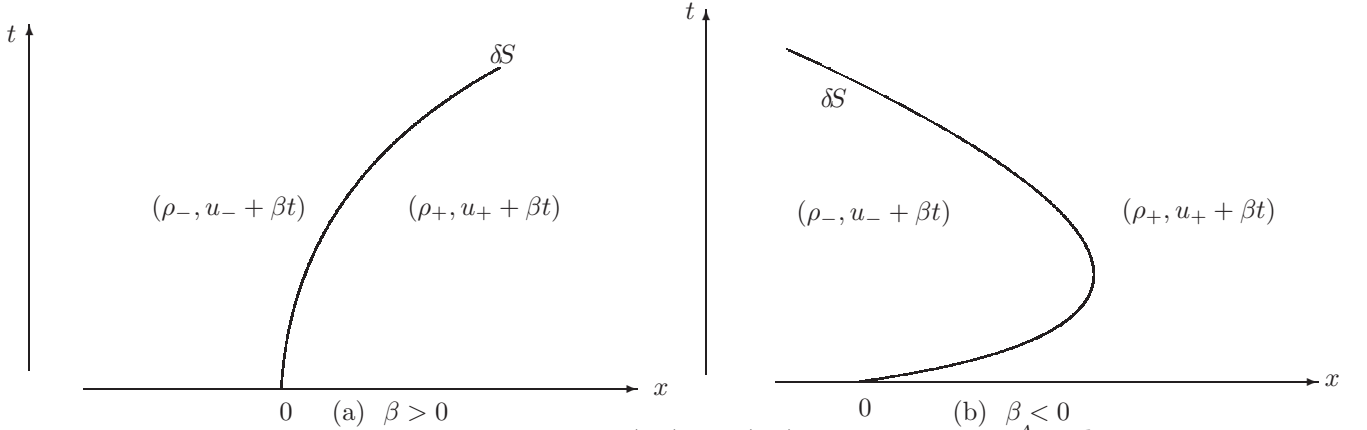
\ \beta<0$} \put(90,5){\makebox(0,0)[cc] {Fig.4.2 The delta shock
wave solution to (1.1) and (1.2) when
$u_+<u_--\frac{A}{\rho_-^\alpha}$ and $v_\delta>0$, }}
\put(90,1){\makebox(0,0)[cc] { where $v_\delta$ is given by (3.28)
for $\rho_-\neq\rho_+$ and (3.30) for $\rho_-= \rho_+$.}}
 \qbezier(115,13.5)(182,22.75)(105.5,57.75)
\end{picture}

For $\beta>0$ (see Fig.4.2(a)), there exists an inverse function of
$x(t)$ globally in the time $t$, which may be written in the form
$$t(x)=\sqrt{\frac{v_{\delta}^{2}}{\beta^{2}}+\frac{2x}{\beta}}-\frac{v_{\delta}}{\beta}.$$
Otherwise, for $\beta<0$ (see Fig.4.2(b)), there is a critical point
$(-\frac{v_{\delta}^{2}}{2\beta},-\frac{v_{\delta}}{\beta})$ on the
delta shock wave curve such that $x'(t)$ change its sign when across
the critical point. Thus, the inverse function of $x(t)$ is needed
to find respectively for $t\leq-\frac{v_{\delta}}{\beta}$ and
$t>-\frac{v_{\delta}}{\beta}$, which enable us to have
\begin{equation}\nonumber
t(x)=\left\{\begin{array}{ll}
-\sqrt{\frac{v_{\delta}^{2}}{\beta^{2}}+\frac{2x}{\beta}}-\frac{v_{\delta}}{\beta},\ \ \ t\leq-\frac{v_{\delta}}{\beta},\\[4pt]
\sqrt{\frac{v_{\delta}^{2}}{\beta^{2}}+\frac{2x}{\beta}}-\frac{v_{\delta}}{\beta},\
\  \ \ \ t>-\frac{v_{\delta}}{\beta}.
\end{array}\right.
\end{equation}

Without loss of generality, let us assume that $\beta>0$ for
simplicity. Actually, the other situation can be dealt with
similarly. Under our assumption, it follows from (\ref{4.11}) that
the position of delta shock wave satisfies $x=x(t)>0$ for all the
time. It follows from (\ref{4.10}) that
\begin{eqnarray}
\frac{d\psi(x(t),t)}{dt}&=&\psi_t(x(t),t)+\frac{dx(t)}{dt}\psi_x(x(t),t)\nonumber\\
&=&\psi_t(x(t),t)+(v_\delta+\beta t)\psi_x(x(t),t)\nonumber\\
&=&\psi_t(x(t),t)+u_\delta(t)\psi_x(x(t),t)\nonumber.
\end{eqnarray}

By exchanging the ordering of integrals and using integration by
parts, we have

\begin{eqnarray}
I&=&\int_{0}^{\infty}\int^{\infty}_{t(x)}\rho_{-}( u_{-}+\beta
t-\frac{A}{\rho_{-}^\alpha})\psi_{t}dt dx
+\int_{0}^{\infty}\int^{\infty}_{t(x)}\rho_{-}(u_{-}+\beta t)(u
_{-}+\beta t-\frac{A}{\rho_{-}^\alpha})\psi_{x}dt dx \nonumber\\
&&+\int_{0}^{\infty}\int_{0}^{t(x)}\rho_{+}( u_{+}+\beta
t-\frac{A}{\rho_{+}^\alpha})\psi_{t}dt
dx+\int_{0}^{\infty}\int_{0}^{t(x)}\rho_{+}(u_{+}+\beta t)(u
_{+}+\beta t-\frac{A}{\rho_{+}^\alpha}) \psi_{x}dtdx\nonumber\\
&&+\int_{0}^{\infty}w_{0}t(v_{\delta}+\beta
t)d\psi(x(t),t)\nonumber\\
&=&\int_{0}^{\infty}(\rho_{+}( u_{+}+\beta
t(x)-\frac{A}{\rho_{+}^\alpha})-\rho_{-}( u_{-}+\beta
t(x)-\frac{A}{\rho_{-}^\alpha}))\psi(x,t(x))dx\nonumber\\
&&+\int_{0}^{\infty}(\rho_{-}(u_{-}+\beta t)( u_{-}+\beta
t-\frac{A}{\rho_{-}^\alpha})-\rho_{+}(u_{+}+\beta t)( u_{+}+\beta
t-\frac{A}{\rho_{+}^\alpha}))\psi(x(t),t)dt \nonumber\\
&&-\int_{0}^{\infty}\int^{\infty}_{t(x)}\beta\rho_{-}\psi dt
dx-\int_{0}^{\infty}\int_{0}^{t(x)}\beta\rho_{+}\psi dt
dx-\int_{0}^{\infty}w_{0}(v_{\delta}+2\beta
t)\psi(x(t),t)dt\nonumber\\
&=&\int_{0}^{\infty}C(t)\psi(x(t),t)dt-\beta(\int_{0}^{\infty}\int_{-\infty}^{x(t)}\rho_{-}\psi
dx dt+\int_{0}^{\infty}\int^{\infty}_{x(t)}\rho_{+}\psi dx
dt)\label{4.12},
\end{eqnarray}
in which \begin{eqnarray}C(t)&=&(\rho_{+}( u_{+}+\beta
t-\frac{A}{\rho_{+}^\alpha})-\rho_{-}( u_{-}+\beta
t-\frac{A}{\rho_{-}^\alpha}))(v_{\delta}+\beta t)\nonumber\\
&&+(\rho_{-}(u_{-}+\beta t)( u_{-}+\beta
t-\frac{A}{\rho_{-}^\alpha})-\rho_{+}(u_{+}+\beta t)( u_{+}+\beta
t-\frac{A}{\rho_{+}^\alpha}))\nonumber\\
&&-w_{0}(v_{\delta}+2\beta t)\nonumber.
\end{eqnarray}
By a tedious calculation, we have
\begin{equation}\label{4.13}
C(t)=-\beta w_{0}t=-\beta w(t).
\end{equation}
Thus, it can be concluded from (\ref{4.12}) and (\ref{4.13})
together that the second equality in (\ref{4.3}) holds in the sense
of distributions. The proof is completed.

\begin{rem}\label{rem:4.1}
%\textbf{ Remark 3.3.}
If both $u_{+}< u_{-}-\frac{A}{\rho_{-}^\alpha}$ and
$\rho_{+}=\rho_{-}$ are satisfied, then the delta shock solution to
the Riemann problem (\ref{1.2}) and (\ref{1.5}) can be expressed in
the form (\ref{4.4}) where
\begin{equation}\label{4.14}
\sigma(t)=u_{\delta}(t)=\frac{1}{2}(u_++u_--\frac{A}{\rho_{-}^\alpha})+\beta
t,\
x(t)=\frac{1}{2}(u_++u_--\frac{A}{\rho_{-}^\alpha})t+\frac{1}{2}\beta
t^{2},\ w(t)=(\rho_-u_--\rho_+ u_+)t.
\end{equation}
The process of proof is completely similar and we omit the details.
\end{rem}

\begin{rem}\label{rem:4.2}   If $u_{+}=u_{-}-\frac{A}{\rho_{-}^\alpha}$, then the
delta shock solution to the Riemann problem (\ref{1.2}) and
(\ref{1.5}) can be also expressed as the form in Theorem
\ref{thm:4.1} and Remark \ref{rem:4.1}. The process of proof is easy
and we omit the details.
\end{rem}

\section{The flux approximation limits of Riemann solutions to (\ref{1.2})}

In this section, we are concerned that the flux approximation limits
of Riemann solutions to (\ref{1.2}) and (\ref{1.5}) converge to the
corresponding ones to (\ref{1.1}) and (\ref{1.5}) or not when the
perturbation parameter $A$ tends to zero. According to the relations
between $u_-$ and $u_+$, we will divide our discussion into the
following three cases:

 (1) $u_-<u_+$; \ \ \ \ \ \     (2) $u_-=u_+$;  \ \ \ \  \ \    (3) $u_->u_+$.

\textbf{Case 5.1.}  $u_-<u_+$

In this case, $(\rho_+,u_+)\in$ I in the $(\rho,v)$ plane, so the
Riemann solutions to (\ref{1.2}) and (\ref{1.5}) $R_1+J$ is given by
(\ref{4.1}), where $x_{1}^-(t),\ x_{1}^+(t)$ and $x_{2}(t)$ are
given by (\ref{3.15}) and (\ref{3.16}) respectively and the states
$(\rho_1,v_1)$ and $(\rho_*,v_*)$ can be calculated as (\ref{3.6})
and (\ref{3.13}). From (\ref{3.6}) and (\ref{3.13}) we have
$$\lim\limits_{A\rightarrow0}\rho_1=\lim\limits_{A\rightarrow0}\big(\frac{A(1-\alpha)}{\frac{x}{t}-\beta
t-w_{-}}\big)^{\frac{1}{\alpha}}=0,$$
$$\lim\limits_{A\rightarrow0}\rho_*=\lim\limits_{A\rightarrow0}\big(\frac{A}{u_+-u_-+\frac{A}{\rho_-^\alpha}}\big)^\frac{1}{\alpha}=0,$$
 which indicate the occurrence
of the vacuum states. Furthermore, the Riemann solutions to
(\ref{1.2}) and (\ref{1.5}) converge to

\begin{equation}\label{5.1}
\lim\limits_{A\rightarrow0}(\rho,u)(x,t)=\left\{\begin{array}{ll}
(\rho_-,u_-+\beta t),\ \ \ \ \ \ \ \ \ \ \ x<u_{-}t+\frac{1}{2}\beta t^{2},\\
vacuum,\ \ \ \ \ \ \ \ \ \ \ \ \ \ \ \ \ u_{-}t+\frac{1}{2}\beta t^{2}<x<u_{+}t+\frac{1}{2}\beta t^{2},\\
(\rho_+,u_++\beta t),\ \ \ \ \ \ \ \ \ \ \ x>u_{+}t+\frac{1}{2}\beta
t^{2},
\end{array}\right.
\end{equation}
which is exactly the corresponding Riemann solutions to the
pressureless Euler equations with the same source term and the same
initial data.

\textbf{Case 5.2.}  $u_-=u_+$

In this case, $(\rho_+,u_+)$ is on the $J$ curve in the $(\rho,v)$
plane, so the Riemann solutions to (\ref{1.2}) and (\ref{1.5}) is
given as follows:
\begin{equation}\label{5.2}
(\rho,u)(x,t)=\left\{\begin{array}{ll}
(\rho_-,u_-+\beta t),\ \ \ \ \ \ \ \ \ \ \ x<u_{-}t+\frac{1}{2}\beta t^{2},\\
(\rho_+,u_++\beta t),\ \ \ \ \ \ \ \ \ \ \ x>u_{+}t+\frac{1}{2}\beta
t^{2},
\end{array}\right.
\end{equation}
which is exactly the corresponding Riemann solutions to the
pressureless Euler equations with the same source term and the same
initial data .

\textbf{Case 5.3.}  $u_->u_+$
\begin{lem}\label{lem:5.1}
If $u_->u_ +$, there exists $A_{1}>A_{0}>0$, such that
$(\rho_+,u_+)\in$ $\rm{I\!I}$ as $A_{0}<A<A_{1}$, and
$(\rho_+,u_+)\in$ $\rm{I\!I\!I}$ as $A\leq A_{0}$.
\end{lem}

\noindent\textbf{Proof.} If $(\rho_+,u_+)\in$ I\!I , then
$0<u_{-}-\frac{A}{\rho_{-}^\alpha}<u_+<u_-$, which gives
$\rho_{-}^\alpha(u_--u_+)<A<\rho_{-}^\alpha u_{-}$. Thus we take
$A_0=\rho_{-}^\alpha(u_--u_+)$ and $A_1=\rho_{-}^\alpha u_-$, then
$(\rho_+,u_+)\in$ I\!I as $A_{0}<A<A_{1}$ and $(\rho_+,u_+)\in$
I\!I\!I as $A\leq A_{0}$.

When $A_{0}<A<A_{1}$, $(\rho_+,u_+)\in$ I\!I in the $(\rho,v)$
plane, so the Riemann solutions to (\ref{1.2}) and (\ref{1.5}) is
given by (\ref{4.2}), where $x_{1}(t)$ and $x_{2}(t)$ are given by
(\ref{3.18}) and (\ref{3.16}) respectively and the states
$(\rho_*,v_*)$ can be calculated as (\ref{3.13}). From (\ref{3.13})
we have
$$\lim\limits_{A\rightarrow A_0}\rho_*=\lim\limits_{A\rightarrow A_0}\big(\frac{A}{u_+-u_-+\frac{A}{\rho_-^\alpha}}\big)^\frac{1}{\alpha}
=\lim\limits_{A\rightarrow
A_0}\big(\frac{\rho_-A}{A-A_{0}}\big)^\frac{1}{\alpha}=\infty.$$

Furthermore, we have the following result.
\begin{lem}\label{lem:5.2}
Let $\frac{dx_1(t)}{dt}=\sigma_1(t)$,
$\frac{dx_2(t)}{dt}=\sigma_2(t)$, then we have
\begin{equation}\label{5.3}
\lim\limits_{A\rightarrow A_0}v_*+\beta t=\lim\limits_{A\rightarrow
A_0}\sigma_1(t)=\lim\limits_{A\rightarrow
A_0}\sigma_2(t)=(u_{-}-\frac{A_0}{\rho_{-}^\alpha})t+\beta
t=u_++\beta t=:\sigma(t),
\end{equation}
\begin{equation}\label{5.4}
\lim\limits_{A\rightarrow A_0}\int_{x_1(t)}^{x_2(t)}\rho_*dx=
\rho_-(u_--u_+)t,
\end{equation}
\begin{equation}\label{5.5}
\lim\limits_{A\rightarrow A_0}\int_{x_1(t)}^{x_2(t)}\rho_*(v_*+\beta
t)dx=\rho_-(u_--u_+)(u_++\beta t) t.
\end{equation}
\end{lem}

\noindent\textbf{Proof.} (\ref{5.3}) is obviously true. We will only
prove (\ref{5.4}) and (\ref{5.5}).

\begin{equation}\nonumber
\lim\limits_{A\rightarrow
A_0}\int_{x_1(t)}^{x_2(t)}\rho_*dx=\lim\limits_{A\rightarrow
A_0}\rho_*(x_2(t)-x_1(t))=\lim\limits_{A\rightarrow
A_0}\rho_*(u_+-\frac{\rho_{\ast}v_{\ast}-\rho_-u_-}{\rho_{\ast}-\rho_-})t=\rho_-(u_--u_+)t,
\end{equation}
\begin{eqnarray}
\lim\limits_{A\rightarrow A_0}\int_{x_1(t)}^{x_2(t)}\rho_*(v_*+\beta
t)dx=(u_++\beta t)\lim\limits_{A\rightarrow
A_0}\int_{x_1(t)}^{x_2(t)}\rho_*dx=\rho_-(u_--u_+)(u_++\beta t)
t\nonumber.
\end{eqnarray}
The proof is completed.

It can be concluded from Lemma \ref{lem:5.2} that the curves of the
shock wave $S_{1}$ and the contact discontinuity $J$ will coincide
when $A$ tends to $A_0$ and the delta shock waves will form.  Next
we will arrange the values which gives the exact position,
propagation speed and strength of the delta shock wave according to
Lemma \ref{lem:5.2}.

 From (\ref{5.4}) and (\ref{5.5}), we let
\begin{equation}\label{5.6}
w(t)=\rho_-(u_--u_+) t,,
\end{equation}
\begin{equation}\label{5.7}
w(t)u_\delta(t)=\rho_-(u_--u_+)(u_++\beta t)t,
\end{equation}
then
\begin{equation}\label{5.8}
u_\delta(t)=(u_++\beta t),
\end{equation}
which is equal to $\sigma(t)$. Furthermore, by letting
$\frac{dx(t)}{dt}=\sigma(t)$, we have
\begin{equation}\label{5.9}
x(t)=u_+ t+\frac{1}{2}\beta t^2.
\end{equation}

From (\ref{5.6})-(\ref{5.9}), we can see that the quantities defined
above are exactly consistent with those given by
(\ref{3.27})-(\ref{3.31}) or (\ref{4.10}) in which we take  $A=A_0$.
Thus, it uniquely determines that the limits of the Riemann
solutions to the system (\ref{1.2}) and (\ref{1.5}) when
$A\rightarrow A_0$ in the case $(\rho_+,u_+)\in$ I\!I is just the
delta shock solution of (\ref{1.2}) and (\ref{1.5}) in the case
$(\rho_+,u_+)\in$ $S_\delta$, where $S_\delta$ is actually the
boundary between the regions I\!I and I\!I\!I. So we get the
following results in the case $u_+<u_-$.

\begin{thm}\label{thm:5.1}
If $u_{+}< u_{-}$, for each fixed $A$ with $A_{0}<A<A_{1}$,
$(\rho_+,u_+)\in$ $\rm{I\!I}$, assuming that $(\rho,u)$ is a
solution containing a shock wave $S_1$ and a contact discontinuity
$J$ of (\ref{1.2}) and (\ref{1.5}) which is constructed in Section
4, it is obtained that when $A\rightarrow A_0$, $(\rho,u)$ converges
to a delta shock wave solution of (\ref{1.2}) and (\ref{1.5}) when
$A=A_0$.
\end{thm}

When $A\leq A_{0}$, $(\rho_+,u_+)\in$ I\!I\!I, so the Riemann
solutions to (\ref{1.2}) and (\ref{1.5}) is given by (\ref{4.4})
with (\ref{4.10}) or (\ref{4.14}), which is a delta shock wave
solution. It is easy to see that when $A\rightarrow 0$, for
$\rho_+\neq \rho_-$,
$$x(t)\rightarrow \sigma_0 t+\frac{1}{2}\beta
t^{2},\ \ w(t)\rightarrow\sqrt{\rho_+\rho_-}(u_-- u_+)t,\ \
\sigma(t)=u_{\delta}(t)\rightarrow  \sigma_0 +\beta t,$$ where
 $\sigma_0=\frac{\sqrt{\rho_-}u_-+\sqrt{\rho_+}u_+}{\sqrt{\rho_-}+\sqrt{\rho_+}}$,
for $\rho_+= \rho_-$,
$$x(t)\rightarrow \frac{1}{2}(u_++u_-)t+\frac{1}{2}\beta
t^{2},\ \ w(t)\rightarrow \rho_+(u_-- u_+)t,\ \
\sigma(t)=u_{\delta}(t)\rightarrow  \frac{1}{2}(u_++u_-)+\beta t,$$
which is exactly the corresponding Riemann solutions to the
pressureless Euler equations with the same source term and the same
initial data \cite{Shen1}. Thus, we have the following result:

\begin{thm}\label{thm:5.2}
If $u_{+}< u_{-}$, for each fixed $A<A_0$, $(\rho_+,u_+)\in$
$\rm{I\!I\!I}$, assuming that $(\rho,u)$ is a delta shock wave
solution of (\ref{1.2}) and (\ref{1.5}) which is constructed in
Section 4, it is obtained that when $A\rightarrow 0$, $(\rho,u)$
converges to a delta shock wave solution to the pressureless Euler
equations with the same source term and the same initial data
\cite{Shen1}.
\end{thm}

Now we summarize the main result in this section as follows.
\begin{thm}\label{thm:5.3}
As the perturbed parameter $A\rightarrow0$, the Riemann solutions to
the approximated nonhomogeneous system (\ref{1.2}) tend to the three
kinds of Riemann solutions to the Riemann solutions to
nonhomogeneous pressureless Euler equations with the same source
term and the same initial data, which include a delta shock wave and
a vacuum state. That is to say, the Riemann solutions to the
transportation equations with Coulomb-like friction is stable under
this kind of flux perturbation.
\end{thm}

\section{Conclusions and Discussions}
It can be seen from the above discussions that the limits of
solutions to the Riemann problem (\ref{1.2}) and (\ref{1.5})
converge to the corresponding ones of the Riemann problem
(\ref{1.1}) and (\ref{1.5}) as $A\rightarrow0$. The approximated
system (\ref{1.2}) is strictly hyperbolic. Although the
characteristic field for $\lambda_{1}$ is genuinely nonlinear, the
characteristic field for $\lambda_{2}$ is still linearly degenerate
and (\ref{1.2}) still belongs to the Temple class. Thus, this
perturbation does not totally change the structure of Riemann
solutions to (\ref{1.1}).

If we also consider the approximation of the flux functions for
(\ref{1.1}) in the form
\begin{equation}\label{6.1}
\left\{\begin{array}{ll}
\rho_t+(\rho u)_x=0,\\
(\rho( u+\frac{1}{1-\alpha}P))_t+(\rho u( u+P))_x=\beta\rho,
 \end{array}\right.
\end{equation}
where $P$ is also given by (\ref{1.3}). We can check that
(\ref{6.1}) has two different eigenvalues $\lambda=u\pm\sqrt{\alpha
B\rho^{-\alpha}u}$, and the characteristic fields for both the two
eigenvalues are genuinely nonlinear. Hence, (\ref{6.1}) is strictly
hyperbolic and by simple calculation, it can be seen that
(\ref{6.1}) does not belong to the Temple class anymore. It is clear
to see that the Riemann solutions for the approximated system
(\ref{6.1}) have completely different structures from those for the
original system (\ref{1.1}). Similar to
\cite{Shen2,Shen-Sun,Sheng-Wang-Yin,Sun,Wang}, we can construct the
Riemann solutions to the Riemann problem (\ref{6.1}) and (\ref{1.5})
in all situations and prove them converge to the corresponding ones
to the Riemann problem (\ref{1.1}) and (\ref{1.5}) as $
A\rightarrow0$.

%
%
%
%
%%% The Appendices part is started with the command \appendix;


%% References


%\bibliographystyle{elsarticle-num}



\end{document}